\documentclass[12pt, a4paper]{amsart}
\usepackage{amssymb}
\usepackage{mathrsfs}

\headsep=1cm \numberwithin{equation}{section}
\newtheorem{thm}{Theorem}[section]
\newtheorem{cor}[thm]{Corollary}
\newtheorem{lem}[thm]{Lemma}
\newtheorem{prop}[thm]{Proposition}
\theoremstyle{definition}
\newtheorem{defn}[thm]{Definition}
\newtheorem{rem}[thm]{Remark}
\numberwithin{equation}{section}

\newcommand{\nat}{\mathbb N}

\newcommand{\Supp}{\rm Supp }
\newcommand{\depth}{\rm depth }
\newcommand{\Ht}{\rm ht }

\newcommand{\Hom}{\rm Hom }
\newcommand{\Ext}{\rm Ext}
\newcommand{\p}{\frak p }
\newcommand{\q}{\frak q }

\newcommand{\fa}{\frak a }
\newcommand{\fb}{\frak b }
\newcommand{\fc}{\frak c }

\newcommand{\Spec}{\rm Spec }

\newcommand{\gdim}{\rm Gdim }
\newcommand{\grade}{\rm grade }

\newcommand{\To}{\longrightarrow}
\newcommand{\V}{\rm Var}

\begin{document}
\title{On the generalization of Faltings' Annihilator Theorem }
\author{Mohammad-Reza Doustimehr and Reza Naghipour$^*$}
\address{Department of Mathematics, University of Tabriz, Tabriz, Iran;
and School of Mathematics, Institute for Research in Fundamental
Sciences (IPM), P.O. Box 19395-5746, Tehran, Iran.}
\email{naghipour@ipm.ir} \email {naghipour@tabrizu.ac.ir}
\email{m\b{ }doustimehr@tabrizu.ac.ir}
\thanks{ 
This research  was been in part supported by a grant from IPM.\\
$^*$Corresponding author: e-mail: {\it naghipour@ipm.ir} (Reza Naghipour)}%

\begin{abstract}
Let $R$ be a commutative Noetherian ring and let $n$ be a non-negative integer.
In this article, by using the theory of Gorenstein dimensions, it is shown that whenever
$R$ is a homomorphic image of a Noetherian Gorenstein ring,  then  the invariants
$\inf\{i\in\nat_0|\, {\dim\Supp}(\fb^tH_{\fa}^i(M))\geq n\text{ for all } t\in\nat_0\}$
and  $\inf\{\lambda_{\fa R_{\p}}^{\fb R_{\p}}(M_{\p})|\,\p\in {\rm Spec} \, R \text{ and } \dim R/ \p\geq
n\}$ are equal, for every finitely generated $R$-module $M$ and for every ideals $\frak a, \frak b$ of $R$
with  $\frak b\subseteq \frak a$. This generalizes the Faltings' Annihilator Theorem [G. Faltings, {\it \"Uber die Annulatoren
lokaler Kohomologiegruppen}, Arch. Math. {\bf30} (1978) 473-476].

\end{abstract}
\maketitle

\section{Introduction}
Throughout this paper, let $R$ denote a commutative Noetherian ring
(with identity) and $\frak a$ an ideal of $R$. For an $R$-module $M$, the
$i$th local cohomology module of $M$ with support in $\V(\fa)$
is defined as:
$$H^i_{\frak a}(M) = \underset{n\geq1} {\varinjlim}\,\, \Ext^i_R(R/\frak a^n, M).$$
Local cohomology was first defined and studied by Grothendieck. We refer the reader to \cite{BS} or \cite{Gr1} for more
details about local cohomology. An important theorem in local
cohomology is Faltings' Annihilator Theorem \cite{Fa1} for local cohomology modules, which
states that, if $R$ is a homomorphic image of a regular ring or $R$ has a dualizing complex, then the invariants
$f_{\fa}^{\fb}(M)$ and $\lambda_{\fa}^{\fb}(M)$ are equal, for every choice of the finitely generated $R$-module $M$ and
for every choice of the ideals $\frak a, \frak b$ of $R$ with  $\frak b\subseteq \frak a$, where
 $f_{\fa}^{\fb}(M)=\inf\{i\in\nat_0\mid\fb\not\subseteq {\rm Rad}(0:H_{\fa}^i(M))\}$ (resp.
$\lambda_{\fa}^{\fb}(M)=\inf\{ {\depth} M_{\p}+ {\Ht}(\fa+\p)/{\p}\mid \p\in{\Spec}R\setminus\V(\fb)\}$)
is the $\frak b$-{\it finiteness dimension of $M$ relative to} $\frak a$ (resp.  the  $\frak b$-{\it minimum  $\frak a$-adjusted depth of $M$}),
see  \cite[Definitions 9.1.5 and 9.2.2]{BS}.

Recently, Khashyarmanesh and Salarian in \cite{KS}, gave a very elegant generalization of  the Faltings' Annihilator Theorem over Gorenstein rings. It is well-known that

\begin{align*}
f_{\fa}^{\fb}(M)&=\inf\{i\in\nat_0\mid \fb^tH_{\fa}^i(M)\not=0 \text{ for all }t\in\nat_0\}\\
&=\inf\{i\in\nat_0\mid {\dim\Supp}\,\fb^tH_{\fa}^i(M)\geq 0 \text{ for
all }t\in\nat_0\}.
\end{align*}
and
\begin{align*}
\lambda_{\fa}^{\fb}(M)&=\inf\{\lambda_{\fa R_{\p}}^{\fb R_{\p}}(M_{\p})|\,\p\in {\rm Spec} \, R\}\\
&=\inf\{\lambda_{\fa R_{\p}}^{\fb R_{\p}}(M_{\p})|\, \dim R/ \p\geq 0\}.
\end{align*}

Now, for a non-negative integer $n$, we define the $n$th $\frak b$-{\it finiteness dimension of $M$ relative to} $\frak a$ (resp.  the $n$th $\frak b$-{\it minimum  $\frak a$-adjusted depth of $M$})  by

\begin{align*}
f_{\fa}^{\fb}(M)_n:=\inf\{i\in\nat_0\mid {\dim\Supp}\,\fb^tH_{\fa}^i(M)\geq n \text{ for
all }t\in\nat_0\}.
\end{align*}
(resp.
\begin{align*}
\lambda_{\fa}^{\fb}(M)_n=:\inf\{\lambda_{\fa R_{\p}}^{\fb R_{\p}}(M_{\p})|\, \dim R/ \p\geq n\}).
\end{align*}
Note that  $f_{\fa}^{\fb}(M)_n$ and $\lambda_{\fa}^{\fb}(M)_n$ are either positive integers or $\infty$ and that
$f_{\fa}^{\fb}(M)_0=f_{\fa}^{\fb}(M)$ and $\lambda_{\fa}^{\fb}(M)_0=\lambda_{\fa}^{\fb}(M)$. So it is rather natural
to ask whether Faltings' Annihilator Theorem, as stated in above, generalizes in the obvious way to the invariants
$f_{\fa}^{\fb}(M)_n$ and $\lambda_{\fa}^{\fb}(M)_n$.  More precisely, as a main result of this paper, we prove the following:
\begin{thm}
Assume that $R$ is a homomorphic image of a Gorenstein ring. Let
$\fa$ and $\fb$ be ideals of $R$ such that $\fb\subseteq\fa$, and
let $M$ be a finitely generated $R$-module. Then, for every non-negative integer $n$,
$$f_{\fa}^{\fb}(M)_n=\lambda_{\fa}^{\fb}(M)_n.$$
\end{thm}
The result in Theorem 1.1 is proved in Theorem 2.14. Our method is based on the theory of Gorenstein dimensions or $G$-dimension.
One of our tools for proving Theorem 1.1 is the following, which will play a key role in this paper.
\begin{prop}
Assume that $R$ is a  Gorenstein ring. Let $\fa$ and $\fb$ be ideals of $R$ such that $\fb\subseteq\fa$, and
let $M$ be a finitely generated $R$-module. Then, for every non-negative integer $n$,
$$f_{\fa}^{\fb}(M)_n=\lambda_{\fa}^{\fb}(M)_n.$$
\end{prop}

Throughout this paper, $R$ will always be a commutative Noetherian
ring with non-zero identity, $\fa, \fb$ will denote ideals of $R$ and $M$ will denote a finitely generated $R$-module.
For any ideal $\frak a$ of $R$, we denote $\{\frak p \in {\rm Spec}\,R:\, \frak p\supseteq \frak a \}$ by $V(\fa)$.
Also,  for any ideal $\frak{b}$ of $R$, {\it the radical of} $\frak{b}$, denoted by ${\rm Rad}(\frak{b})$, is defined to
be the set $\{x\in R \,: \, x^n \in \frak{b}$ for some $n \in \mathbb{N}\}$.
For any unexplained notation and terminology concerning Gorenstein dimensions  we refer the reader to \cite{Ch}.

\section{Faltings' Annihilator Theorem}
The main goal of this section is to provide a generalization of the Faltings' Theorem for the annihilation of local cohomology modules over a Gorenstein
ring. The main results are Theorems  2.10 and 2.14.  Firstly, let us, introduce  the concept of Gorenstein dimension. The notion of Gorenstein dimension ($G$-dimension) was introduced by Auslander \cite{A} and was deeply studied by him and Bridger \cite{AB}.  For an $R$-module $L$ the biduality map is the
canonical $R$-homomorphism
$$\delta_L: L\longrightarrow {\rm Hom}_R({\rm Hom}_R(L, R), R),$$
defined  by $\delta_L(x)(f)=f(x)$ for all $f\in {\rm Hom}_R(L, R)$ and $x\in L$.

\begin{defn}\label{def1}(cf. \cite[1.1.2]{Ch})
A finitely generated $R$-module $M$ belongs to the $G$-{\it class} $G(R)$ if and only if
\begin{enumerate}
\item
${\Ext}_R^i(M,R)=0$ for all $i>0$;
\item
${\Ext}_R^i({\Hom}_R(M, R),R)=0$ for all $i>0$;
\item
The biduality map $\delta_M: M\To{\Hom}_R({\Hom}_R(M,R),R)$ is an isomorphism.
\end{enumerate}
\end{defn}
\begin{defn}\label{def2}(cf. \cite[1.2.1]{Ch})
A $G$-{\it resolution} of a finitely generated $R$-module $M$ is a
sequence of modules in $G(R)$,
$$\cdots\To G_l\To G_{l-1}\To\cdots\To G_2\To G_1\To G_0\To 0,$$
which is exact at $G_l$ for all $l>0$ and has $G_0/{\rm Im}(G_1\To
G_0)\simeq M$. That is, there is an exact sequence
$$\cdots\To G_l\To G_{l-1}\To\cdots\To G_2\To G_1\To G_0\To M\To 0.$$
The resolution is said to be of finite length $n$ if $G_n\not=0$ and
$G_l=0$ for all $l>n$.
\end{defn}

Note that every finitely generated $R$-module has a resolution by
finitely generated free modules and, thereby, a $G$-resolution.

\begin{defn}\label{def3}(cf. \cite[1.2.3]{Ch})
A finitely generated $R$-module $M$ is said to have {\it finite
$G$-dimension}, if it has a $G$-resolution of finite length. We set
${\gdim} 0=-\infty$ and for $M\not=0$, we define $G$-dimension of
$M$ as follows: for any positive integer $n$, we say that $M$ has
$G$-dimension at most $n$, and write ${\gdim}_R M\leqslant n$ if and
only if $M$ has a $G$-resolution of length n. If $M$ has no
$G$-resolution of finite length, then we say that it has infinite
$G$-dimension and write ${\gdim}_R M=\infty$.
\end{defn}
\begin{cor}\label{cor4}
 A commutative Noetherian local ring $R$ is
Gorenstein if and only if every finitely generated $R$-module has
finite $G$-resolution.
\end{cor}

\proof See \cite[Corollary 2]{Go}.\qed
\begin{rem}\label{rem5}
 Let $M$ be a finitely generated $R$-module. For each $t\in\nat_0\cup\{-\infty\}$, put
$C_t(M)=\{\p\in{\Spec}\,R\mid {\gdim} M_{\p}>t\}$  and $c_t(M)=\bigcap_{\p\in C_t(M)}\p$. Then $C_t(M)$ is a closed subset of
${\Spec}\, R$ (in the Zariski topology) and
$$\sqrt{0:M}=c_{-\infty}(M)\subseteq c_0(M)\subseteq
c_1(M)\subseteq\cdots \subseteq c_t(M)\subseteq \cdots. $$
\end{rem}

We now state and prove some preliminary lemmas and a proposition which help us to conclude the main results.
\begin{lem}\label{lem6}
Let $M$ be a finitely generated $R$-module, and let
$\p\in{\Spec} R$ be such that ${\gdim}M_{\p}<\infty$. Then there
exists $s\in
 R\setminus\p$ such that, for every proper ideal $\fa$ of $R$, we
 have
 $$sH_{\fa}^i(M)=0\text{ \hspace{5mm} for all  } i<{\grade} (\fa, R)-{\gdim}M_{\p}.$$
\end{lem}

\proof Set $h:={\gdim}M_{\p}$. We use induction on $h$. If $h=-\infty$, then $M_{\p}=0$ and the result is clear from
\cite[Lemma 9.4.1]{BS}. When $h=0$, the desired result follows from
\cite[Lemma 2.9]{KS}. We therefore assume, inductively,  that
$h>0$ and the result  has been proved for smaller values of $h$.
There is a non-zero, finitely generated free $R$-module $F$ and an
exact sequence $$0\To N\To F\To M\To 0$$ of $R$-modules and $R$-homomorphisms.  Localization
yields an exact sequence $$0\To N_{\p}\To F_{\p}\To M_{\p}\To 0.$$
Therefore, in view of \cite[Corollary 1.2.9(c)]{Ch}, ${\gdim}N_{\p}=h-1$  and
so, by inductive hypothesis, there exists $s\in R\backslash\p$ such
that, for every proper ideal $\fa$ of $R$, we have
$$sH_{\fa}^i(N)=0\text{ \hspace{5mm} for all }
i<{\grade}(\fa,R)-h+1.$$ Thus $sH_{\fa}^{i+1}(N)=0$ for all
$i<{\grade}(\fa,R)-h$. Let $\fa$ be a proper ideal of $R$ and let
$i\in\nat_0$ with $i<{\grade}(\fa,R)-h$. Now, in view of  the exact sequence
$$H_{\fa}^i(F)\To H_{\fa}^i(M)\To H_{\fa}^{i+1}(N)$$ and \cite[Lemma
6.2.7]{BS}, we have $sH_{\fa}^i(M)=0$. This completes the inductive
step.\qed\\

\begin{lem}\label{lem7}
Let $M$ be a finitely generated $R$-module, and
$t\in\nat_0\cup\{-\infty\}$. Then there exists $n\in\nat$ such that, for
every proper ideal $\fa$ of $R$, we have
$$c_t(M)^nH_{\fa}^i(M)=0\text{ \hspace{5mm} for all } i<{\grade}(\fa,R)-t.$$
\end{lem}

\proof Let $\p\in U:={\Spec}\, R\setminus C_t(M)$. Thus ${\gdim}\,
M_{\p}\leqslant t$. By Lemma 2.6, there exists $s_{\p}\in R\setminus\p$
such that, for every proper ideal $\fa$ of $R$, we have
 $sH_{\fa}^i(M)=0$  for all  $ i<{\grade}(\fa,R)-{\gdim}\,M_{\p}$. Set $\frak g:=\sum_{\p\in U} s_{\p}R$, and
 observe that, for every proper ideal $\fa$ of $R$, we have
 ${\frak g}H_{\fa}^i(M)=0$  for all  $ i<{\grade}(\fa,R)-t$. As
 $c_t(M)\subseteq\sqrt{\frak g}$,  there exists $n\in\nat$
 such that $c_t(M)^n\subseteq{\frak g}$, and the result  now follows from  this.\qed \\

\begin{prop}\label{prop8}
Let $R$ be a Gorenstein  ring,  let $M$ be  a finitely generated
$R$-module, and let  $\fa, \fb$ be  ideals of $R$ such that $\fb \subseteq \fa$. Then, for all
$\q\in V(\fa)$ with $\dim R/\q\geq n$, $$\fb R_{\q}\subseteq
c_{{\Ht}\, {\q}-\lambda_{\fa}^{\fb}(M)_n}(M_{\q}).$$
\end{prop}

\proof Let $\q\in{\V}(\fa)$ with $\dim R/\q\geq n$ and let $\p$ be an arbitrary prime ideal of $R$ such that
$\p R_{\q}\in C_{{\Ht}\, {\q}-\lambda_{\fa}^{\fb}(M)_n}(M_{\q})$. It is enough to show that $\fb
R_{\q}\subseteq \p R_{\q}$. To achieve this, suppose that the contrary is true, i.e., $\fb R_{\q}\not\subseteq \p R_{\q}$,
and look for a contradiction. Then it follows from \cite[Remarks 9.2.3]{BS} that
$$\infty> {\depth}\, M_{\p} + {\Ht}\, (\q R_{\q}+\p R_{\q})/\p
R_{\q}\geq\lambda_{\q R_{\q}}^{\fb R_{\q}}(M_{\q})\geq\lambda_{\fa
R_{\q}}^{\fb R_{\q}}(M_{\q})\geq \lambda_{\fa}^{\fb}(M)_n.$$  Next, since $R_{\q}$
is catenary, it yields that
 $${\Ht}\,\q -{\Ht}\,\p= {\Ht}\,\q R_{\q} -{\Ht}\, \p R_{\q} ={\Ht}\, \q R_{\q}/\p R_{\q}={\Ht}\, (\q R_{\q}+ \p R_{\q})/\p R_{\q}.$$
Also, as  $R_{\p}$ is Gorenstein,  it follows from Auslander-Birdger formula (see \cite[Theorem 1.4.8]{Ch} and Corollary \ref{cor4}, that
$${\depth}\, M_{\p} +{\gdim}\,M_{\p}={\depth}\,R_{\p}={\Ht}\,\p.$$
Consequently $${\Ht}\,\p -{\gdim}\, M_{\p} + {\Ht}\,\q -{\Ht}\,\p\geq   \lambda_{\fa}^{\fb}(M)_n,$$
and so  $${\gdim}\,(M_{\q})_{\p R_{\q}}={\gdim}\, M_{\p}\leqslant {\Ht}\,\q - \lambda_{\fa}^{\fb}(M)_n.$$ Therefore $\p
R_{\q}\not\in C_{{\Ht}\, {\q}-\lambda_{\fa}^{\fb}(M)_n}(M_{\q})$, which is a contradiction. \qed\\

\begin{lem}\label{lem9}
Let $L\To M\To N$ be an exact sequence of $R$-homomorphisms and $R$-modules. Suppose that $n,t$ and $s$  be non-negative integers such that
${\dim\Supp}({\fb}^tL)<n$ and ${\dim\Supp}({\fb}^{s} N)<n$. Then there exists a non-negative integer $l$ such that
${\dim\Supp}({\fb}^lM)<n$.
\end{lem}

\proof Set $l:=t+s$. It is enough to show that for each $\p\in {\Spec}\, R$ with $\dim R/\p\geq n$, we have
 $$(\fb^lM)_{\p}=(\fb R_{\p})^lM_{\p}=0.$$
To do this, let $m\in M_{\p}$ and we consider the exact sequence
$$L_{\p}\stackrel{f}\To M_{\p}\stackrel{g}\To N_{\p},$$ of
$R_{\p}$-modules and $R_{\p}$-homomorphisms.  Now, for each $u\in\fb^{s} R_{\p}$ we have
$g(um)=ug(m)=0$, and so $um\in {\rm Ker} g={\rm Im} f$. Thus there
exists $v\in L_{\p}$ such that $um=f(v)$. Also, for each $w\in\fb^t
R_{\p}$, we have $$wum=wf(v)=f(wv)=0.$$ Hence  $\fb^l R_{\p}m=(\fb^t
R_{\p})(\fb^{s} R_{\p})m=0$, and therefore  $(\fb R_{\p})^lM_{\p}=0$, as required. \qed \\


 We are now ready to state and prove the generalization of the Faltings' Theorem for the annihilation of local cohomology modules over a Gorenstein
ring.
\begin{thm}\label{thm 10} {\rm (Faltings' Annihilator Theorem)}
Assume that $R$ is a Gorenstein ring and $M$ a finitely generated $R$-module.  Let $\fa$ and $\fb$ be ideals
of $R$ such that $\fb\subseteq\fa$. Then, for every non-negative integer $n$,
$$f_{\fa}^{\fb}(M)_n=\lambda_{\fa}^{\fb}(M)_n.$$
\end{thm}

\proof  Let $f_{\fa}^{\fb}(M)_n=h$. Then there exists a non-negative integer $t$  such that for all integers  $i<h$ we have,
 $${\dim\Supp}({\fb}^tH_{\fa}^i(M))<n.$$
Thus  $({\fb}^tH_{\fa}^i(M))_{\p}=0$ for all $i<h$ and for all $\p\in {\Spec}\, R$ with $\dim R/\p\geq n$, and so $f_{\fa
R_{\p}}^{\fb R_{\p}}(M_{\p})\geq h$ for all $\p\in{\Spec}\,R$ with $\dim R/\p\geq n$. Hence, on use of \cite[Theorem 9.3.5]{BS}, we have
\begin{align*}
f_{\fa}^{\fb}(M)_n=h&\leqslant\inf\{f_{\fa R_{\p}}^{\fb
R_{\p}}(M_{\p})\mid \p\in{\Spec}\, R\text{ and }\dim R/\p\geq
n\}\\&\leqslant\inf\{\lambda_{\fa R_{\p}}^{\fb R_{\p}}(M_{\p})\mid
\p\in{\Spec}\, R \text{ and }\dim R/\p\geq n\}=\lambda_{\fa}^{\fb}(M)_n.
\end{align*}

Now, for proving the inequality  $f_{\fa}^{\fb}(M)_n\geq\lambda_{\fa}^{\fb}(M)_n$, there are two
cases to consider.\\

{\bf Case 1.} Suppose that $\dim R/\fa\leq n$. Then the set
$$\mathcal{T}:=\{\p\,|\,\, \p\in V(\fa)\,  \text{and}\,  \dim R/\p=n\},$$ is  finite. Let
$\mathcal{T}=\{\p_1,\dots, \p_h\}$ and consider an integer $j$ with $1\leqslant
j\leqslant h$. Set $t_j={\Ht}\,\p_j-\lambda_{\fa}^{\fb}(M)_n$. Since
${\Ht}\,\p_j={\grade}(\p_jR_{\p_j},R_{\p_j})$,
 we can deduce from Lemma  \ref{lem7} that there exists $s_j\in\nat$ such that
 $$c_{t_j}(M_{\p_j})^{s_j}H_{\p_jR_{\p_j}}^i(M_{\p_j})=0 \text{\hspace{15mm} for
 all } i<{\Ht}\,{\p_j}-t_j=\lambda_{\fa}^{\fb}(M)_n.$$
Moreover,  in view of Proposition 2.8, $\fb R_{\p_j}\subseteq c_{t_j}(M_{\p_j})$. Set
$t=\max\{s_1,\dots, s_h\}$. Then
$$({\fb}^tH_{\fa}^i(M))_{\p_j}=(({\fb R_{\p_j}})^tH_{\fa
R_{\p_j}}^i(M_{\p_j}))\subseteq
c_{t_j}(M_{\p_j})^tH_{\p_jR_{\p_j}}^i(M_{\p_j})=0 \text{\hspace{5mm}
for all }  i<\lambda_{\fa}^{\fb}(M)_n.$$
Hence $\p_1,\dots, \p_h\not\in{\Supp}({\fb}^tH_{\fa}^i(M))$,  and so
${\dim\Supp}({\fb}^tH_{\fa}^i(M))\leqslant n-1$. Therefore
$f_{\fa}^{\fb}(M)_n\geq \lambda_{\fa}^{\fb}(M)_n$, as required.\\

{\bf Case 2.} Now, suppose that $\dim R/\fa>n$ and we show that $f_{\fa}^{\fb}(M)_n\geq \lambda_{\fa}^{\fb}(M)_n$.
Suppose, the contrary, that $f_{\fa}^{\fb}(M)_n<\lambda_{\fa}^{\fb}(M)_n$
 and look for a contradiction. To this end, as $R$ is Noetherian,  we can (and do) assume that $\fa$ is a maximal element of the set
$$\Sigma:=\{\fa'\in \mathscr{I}(R)\mid\text{ there exists  $\fb'\in \mathscr{I}(R)$ with $\fb'\subseteq\fa'$ and }
f_{\fa'}^{\fb'}(M)_n<\lambda_{\fa'}^{\fb'}(M)_n \},$$
where $\mathscr{I}(R)$ denotes  the set of all ideals of $R$.
 Let $\p_1,\dots, \p_h$ be the distinct minimal primes of $\fa$. Then
 there is an integer $i$ with $1\leq i\leq h$ such that
 $\dim R/\fa=\dim R/\p_i$. By \cite[Exercise 9.4.9]{BS} and Proposition \ref{prop8}, we have
$$\fb R_{\p_i}\subseteq c_{({\Ht}\, {\p_i}-\lambda_{\fa}^{\fb}(M)_n)}(M_{\p_i})=
(c_{({\Ht}\, {\p_i}-\lambda_{\fa}^{\fb}(M)_n)}(M))_{\p_i}.$$ Therefore there
exists $u\in R\setminus\p_i$ such that $\fb R_u\subseteq(c_{({\Ht}\,
{\p_i}-\lambda_{\fa}^{\fb}(M)_n)}(M))_u$. As $\dim R/\p_i>0$,
there exists $\q\in{\Spec}\, R$ such that $\p_i\subsetneqq\q$. Let
$v\in\q\setminus\p_i$ and $w\in\cap_{j\not=i}\p_j\setminus\p_i$. Set
$s:=uvw$.  Then in view of  \cite[Exercise 9.4.9]{BS},  we have
$$\fb R_s\subseteq(c_{({\Ht}\, {\p_i}-\lambda_{\fa}^{\fb}(M)_n)}(M))_s=c_{({\Ht}\, {\p_i} R_s-\lambda_{\fa}^{\fb}(M)_n)}(M_s)=c_{({\Ht}\,\fa
R_s-\lambda_{\fa}^{\fb}(M)_n)}(M_s).$$

Since ${\grade}(\fa R_s,R_s)={\Ht}\,\fa R_s$, it follows from Lemma
\ref{lem7} that there exists $t_1\in\nat$ such that
$$(c_{({\Ht}\,\fa R_s-\lambda_{\fa}^{\fb}(M)_n)}(M_s))^{t_1}H_{\fa
R_s}^i(M_s)=0\text{  \,\, for all $i<\lambda_{\fa}^{\fb}(M)_n$}.$$ Hence
$\fb^{t_1}H_{\fa R_s}^i(M_s)=0$ for all $i<\lambda_{\fa}^{\fb}(M)_n$. Therefore, by virtue of  \cite[Theorem 4.2.1]{BS}, we have
$\fb^{t_1}H_{\fa}^i(M_s)=0$ for all $i<\lambda_{\fa}^{\fb}(M)_n$. Since
$\fa\subsetneqq\fa+Rs\subsetneqq R$, it follows by the `maximality' assumption
on $\fa$ that $\lambda_{\fa+Rs}^{\fb}(M)_n\leq f_{\fa+Rs}^{\fb}(M)_n$. Now, as in view of  \cite[Remark 9.2.3]{BS},
$\lambda_{\fa}^{\fb}(M)_n\leq\lambda_{\fa+Rs}^{\fb}(M)_n$,  it follows that there
exists $t_2\in\nat$ such that
${\dim\Supp} (\fb^{t_2}H_{\fa+Rs}^i(M))\leq n-1$ for all
$i<\lambda_{\fa}^{\fb}(M)_n$. By the exact sequence
$$\cdots\To H_{\fa+Rs}^i(M)\To H_{\fa}^i(M)\To H_{\fa}^i(M_s)\To\cdots$$
and Lemma \ref{lem9}, it follows that there exists $t\in\nat$ such that
$\dim{\rm Supp} (\frak b^tH_{\frak a}^i(M))\leq n-1$ for all
$i<\lambda_{\fa}^{\fb}(M)_n$. Therefore $\lambda_{\fa}^{\fb}(M)_n\leq f_{\fa}^{\fb}(M)_n$, and this contradiction completes the proof.
\qed \\
\begin{lem}\label{lem11}
Let $f: R\To R'$ be a surjective homomorphism of Noetherian rings and let $f^*:{\Spec}\, R'\longrightarrow {\Spec}\, R$ be
the induced map. Then for any $R'$-module $L$, $${\Supp}_R L=f^*({\Supp}_{R'} L).$$
\end{lem}
\proof Follows easily from \cite[Proposition 3.2]{Me} and \cite[Ch. IV, Sec. 3, Proposition 7]{Bo}.\qed \\
\begin{lem}\label{lem12}
Let $f: R\To R'$ be a surjective homomorphism of Noetherian rings, and $\fa$ and $\fb$ be ideals of $R$ such that $\fb\subseteq\fa$.
Let $M'$ be a finitely generated $R'$-module. Then $$f_{\fa}^{\fb}(M')_n=f_{\fa R'}^{\fb R'}( M')_n.$$
\end{lem}

\proof
Let $f^*: {\Spec}\, R'\longrightarrow {\Spec}\, R$ be the induced map.
Then, for each $t\in\nat_0$, we have
$$f^*({\Supp}_{R'}(\fb R')^t H_{\fa R'}^i(M'))= {\Supp}_R(\fb R')^t H_{\fa
R'}^i(M')={\Supp}_R\fb^t H_{\fa}^i(M'),$$ and so
$\dim_{R'}(\fb R')^t H_{\fa R'}^i(M')=\dim_R \fb^t
H_{\fa}^i(M')$. Thus $$f_{\fa}^{\fb}(M')_n=f_{\fa R'}^{\fb R'}( M')_n,$$  as required. \qed \\

Before we state Theorem 2.14 which is the main result of this paper, we give a couple of lemmas that in the 
proof of Theorem 2.14.\\
\begin{lem}\label{lem13}
Let $\fa$ and $\fb$ be ideals of $R$ such that $\fb\subseteq\fa$. Let $M$ be a finitely generated $R$-module, and let $\fc$ be an
ideal of $R$ such that $\fc\subseteq (0:M)$. Then
$$\lambda_{\fa}^{\fb}(M)_n=\lambda_{(\fa+\fc)/\fc}^{(\fb+\fc)/\fc}(M)_n.$$
\end{lem}

\proof In view of  \cite[Lemma 9.2.6]{BS} we have
\begin{align*}
\lambda_{\fa}^{\fb}(M)_n&=\inf\{\lambda_{\fa R_{\p}}^{\fb
R_{\p}}(M_{\p})\mid \p\in{\Spec}\, R\text{ and }\dim R/\p\geq n\} \\
&=\inf\{\lambda_{(\fa R_{\p}+\fc R_{\p})/\fc R_{\p}}^{(\fb R_{\p}+\fc R_{\p})/\fc R_{\p}}(M_{\p})\mid
\p\in{\Spec}\,R\text{ and }\dim R/\p\geq n\} \\
&=\inf\{\lambda_{((\fa+\fc)/\fc)_{\p/\fc}}^{((\fb+\fc)/\fc)_{\p/\fc}}(M_{\p/\fc})\mid
\p/\fc\in{\Spec}\,R/\fc \text{ and }\dim (R/\fc)/(\p/\fc) \geq n\} \\
&=\lambda_{(\fa+\fc)/\fc}^{(\fb+\fc)/\fc}(M)_n, 
\end{align*}
as required.   \qed \\

We are now ready to state and prove the main result of this paper which is a generalization of the Faltings' Theorem for the annihilation of local cohomology modules  whenever the ring $R$ is a homomorphic image of a Noetherian Gorenstein ring.

\begin{thm}\label{thm14}
Assume that $R$ is a homomorphic image of a Gorenstein ring. Let
$\fa$ and $\fb$ be ideals of $R$ such that $\fb\subseteq\fa$, and
let $M$ be a finitely generated $R$-module. Then, for each $n\in\nat_0$,
$$f_{\fa}^{\fb}(M)_n=\lambda_{\fa}^{\fb}(M)_n.$$
\end{thm}

\proof By assumption there is a Gorenstein ring $R'$ and a surjective homomorphism of Noetherian rings $f:R'\To R$. Let $\fa'$
and $\fb'$ be ideals of $R'$ such that $\fa=\fa' R$ and $\fb=\fb'
R$. Then by Lemmas \ref{lem12}, \ref{lem13}, and Theorem 2.10,
$$f_{\fa}^{\fb}(M)_n=f_{\fa' R}^{\fb' R}(M)_n=f_{\fa'}^{\fb'}(M)_n=\lambda_{\fa'}^{\fb'}(M)_n=\lambda_{\fa}^{\fb}(M)_n,$$
as required. \qed \\

\begin{center}
{\bf Acknowledgments}
\end{center}
The authors would like to thank Professor Hossein Zakeri for his reading of the first
draft and valuable discussions. Also, we would like to thank from School of Mathematics, Institute for Research in Fundamental
Sciences (IPM), for its financial support.


\end{document}